\documentclass[12pt]{article}
\usepackage{amsfonts}
\usepackage{a4wide}
\usepackage{amsmath,amscd,amsthm,a4,amssymb}

\begin{document}

\begin{center}
{\Large \textbf{A Chebyshev type alternation theorem for best approximation
by a sum of two algebras}}

\

\textbf{Aida Kh. Asgarova}$^{1}$, \textbf{Ali A. Huseynli}$^{2}$ \textbf{and
Vugar E. Ismailov}$^{3}$ \vspace{1mm}

$^{1,2,3}${Institute of Mathematics and Mechanics, Baku, Azerbaijan} \vspace{%
1mm}

$^{2}${Khazar University, Baku, Azerbaijan} \vspace{1mm}

e-mail: $^{1}$aidaasgarova@gmail.com, $^{2}$alihuseynli@gmail.com, $^{3}$%
vugaris@mail.ru
\end{center}

\smallskip

\textbf{Abstract.} Let $X$ be a compact metric space, $C(X)$ be the space of
continuous real-valued functions on $X$, and $A_{1},A_{2}$ be two closed
subalgebras of $C(X)$ containing constant functions. We consider the problem
of approximation of a function $f\in C(X)$ by elements from $A_{1}+A_{2}$.
We prove a Chebyshev type alternation theorem for a function $u_{0}\in
A_{1}+A_{2}$ to be a best approximation to $f$.

\smallskip



\bigskip

\bigskip

\begin{center}
{\large \textbf{1. Introduction}}
\end{center}

The classical Chebyshev alternation theorem gives a criterion for a
polynomial $P$ of degree not greater than $n$ to be the best uniform
approximation to a continuous real valued function $f$, using the
oscillating nature of the difference $f-P$. More precisely, the theorem
asserts that $P$ is the best uniform approximation to $f$ on $[0,1]$ if and
only if there exist $n+2$ points $t_{i}$ in $[0,1]$ such that

\begin{equation*}
f(t_{k})-P(t_{k})=(-1)^{k}\max_{t\in \lbrack 0,1]}\left\vert
f(t)-P(t)\right\vert ,\text{ }k=1,...,n+2.
\end{equation*}
See the monograph of Natanson \cite{Nat} for a comprehensive commentary on
this theorem. Several general alternation theorems applying to an arbitrary
finite dimensional subspace $M$ of $C(X)$ for $X$ a cell in $\mathbb{R}^{d}$%
, may be found in Buck \cite{Buck}. For the history and various variants of
the Chebyshev alternation theorem consult \cite{Bro}.

In this paper, we prove a Chebyshev type alternation theorem for a best
approximation of a continuous function, defined on a compact metric space,
by sums of two algebras. To make the problem more precise, assume $X$ is a
compact metric space, $C(X)$ is the space of real-valued continuous
functions on $X$, $A_{1}$ and $A_{2}$ are closed subalgebras of $C(X)$
containing constants. For a given function $f\in C(X)$ consider the
approximation of $f$ by elements of $A_{1}+A_{2}$. We ask and answer the
following question: which conditions imposed on $u_{0}\in A_{1}+A_{2}$ are
necessary and sufficient for the equality

\begin{equation*}
\left\Vert f-u_{{0}}\right\Vert =\inf_{u\in A_{1}+A_{2}}\left\Vert
f-u\right\Vert ?\eqno(1.1)
\end{equation*}%
Here $\left\Vert \cdot \right\Vert $ denotes the standard uniform norm in $%
C(X)$. Recall that a function $u_{0}$ satisfying (1.1) is called a best
approximation to $f.$

It should be remarked that approximation problems concerning sums of
algebras were studied in many papers (see Khavinson's monograph \cite{12}
for an extensive discussion). The history of this subject goes back to 1937
and 1948 papers by M.H. Stone \cite{St1,St2}. He considered the most
particular case of the approximation by sums of algebras, namely the case
when only one algebra is involved. A version of the corresponding famous
result, known as the Stone-Weierstrass theorem, states that a subalgebra $%
A\subset C(X)$, which contains a nonzero constant function, is dense in the
whole space $C(X)$ if and only if $A$ separates points of $X$ (that is, for
any two different points $x$ and $y$ in $X$ there exists a function $g\in A$
with $g(x)\neq g(y)$). Density of the sum of two subalgebras $A_{1}$ and $%
A_{2}$ in $C(X)$ (for a compact Hausdorff $X$) was extensively studied in
Marshall and O'Farrell \cite{Mar1,Mar2}. In \cite{Mar1}, they gave a
complete description of measures on $X$ orthogonal to the sum $A_{1}+A_{2}$.
From this description they obtained a geometrical condition which is
equivalent to the density of $A_{1}+A_{2}$ in $C(X)$. The paper \cite{Mar1}
also indicates main difficulties with the sum of more than two algebras.

This paper exploits the same mathematical and geometrically explicit objects
from Marshall and O'Farrell \cite{Mar1,Mar2} for characterization of a best
approximation by a sum of two algebras. To prove our main result we use
various results and ideas of Functional Analysis and General Topology.

Note that the algebras $A_{i}$, in particular cases, turn into algebras of
univariate functions, ridge functions and radial functions. The literature
abounds with the use of ridge functions and radial functions. Ridge
functions and radial functions are defined as multivariate functions of the
\ form $g(\mathbf{a}\cdot \mathbf{x})$ and $g(\left\vert \mathbf{x}-\mathbf{a%
}\right\vert _{e})$ respectively, where $\mathbf{a}\in \mathbb{R}^{d}$ is a
fixed vector, $\mathbf{x}\in \mathbb{R}^{d}$ is the variable, $\mathbf{a}%
\cdot \mathbf{x}$ is the usual inner product, $\left\vert \mathbf{x-a}%
\right\vert _{e}$ is the Euclidean distance between $\mathbf{x}$ and $%
\mathbf{a}$, and $g$ is a univariate function.

\bigskip

\bigskip

\begin{center}
{\large \textbf{2. The main result}}
\end{center}

Let $X$ be a compact metric space, $C(X)$ be the space of real-valued
continuous functions on $X$ and $A_{1}\subset C(X),$ $A_{2}\subset C(X)$ be
two closed algebras that contain the constants. Define the equivalence
relation $R_{i},$ $i=1,2,$ for elements in $X$ by setting

\begin{equation*}
a\overset{R_{i}}{\sim }b\text{ if }f(a)=f(b)\text{ for all }f\in A_{i}.
\end{equation*}%
Then, for each $i=1,2,$ the quotient space $X_{i}=X/R_{i}$ with respect to
the relation $R_{i}$, equipped with the quotient space topology, is compact
and the natural projections $s:X\rightarrow X_{1}$ and $p:X\rightarrow X_{2}$
are continuous. Note that the quotient spaces $X_{1}$ and $X_{2}$ are not
only compact but also Hausdorff (see, e.g., \cite[p.54]{12}). In view of the
Stone-Weierstrass theorem, the algebras $A_{1}$ and $A_{2}$ have the
following set representations
\begin{eqnarray*}
A_{1} &=&\{g(s(x)):~g\in C(X_{1})\}, \\
A_{2} &=&\{h(p(x)):~h\in C(X_{2})\}.
\end{eqnarray*}

We proceed with the definition of \textit{lightning bolts} with respect to
two algebras. These objects are essential for our further analysis.

\bigskip

\textbf{Definition 2.1.} (see \cite{Mar1}) \textit{A finite or infinite
ordered set $l=\{x_{1},x_{2},...\}\subset X$, where $x_{i}\neq x_{i+1}$,
with either $s(x_{1})=s(x_{2}),p(x_{2})=p(x_{3}),s(x_{3})=s(x_{4}),...$ or $%
p(x_{1})=p(x_{2}),s(x_{2})=s(x_{3}),p(x_{3})=p(x_{4}),...$ is called a
lightning bolt with respect to the algebras $A_{1}$ and $A_{2}$.}

\bigskip

In the sequel, we will simply use the term \textquotedblleft
bolt\textquotedblright\ instead of the expression \textquotedblleft
lightning bolt with respect to the algebras $A_{1}$ and $A_{2}$%
\textquotedblright . If in a finite bolt $\{x_{1},...,x_{n},x_{n+1}\}$, $%
x_{n+1}=x_{1}$ and $n$ is an even number, then the bolt $\{x_{1},...,x_{n}\}$
is said to be closed.

Bolts, in the special case when $X\subset \mathbb{R}^{2}$, and $%
A_{1}=\{g(x)\}$ and $A_{2}=\{h(y)\}$, are geometrically explicit objects. In
this case, a bolt is an ordered set $\{x_{1},x_{2},...\}$ in $\mathbb{R}^{2}$
with the line segments $[x_{i},x_{i+1}],$ $i=1,...,n,$ perpendicular
alternatively to the $x$ and $y$ axes. Bolts, in this particular and
simplest case, were first introduced by Diliberto and Straus in \cite{4}.
They were further used in many works devoted to the approximation of
multivariate functions by sums of univariate functions (see \cite{12}).
Bolts appeared in a number of papers with several different names such as
\textit{permissible lines} (see \cite{4}), \textit{paths} (see, e.g., \cite%
{13}), \textit{trips}\ (see, e.g., \cite{Mar2}), \textit{links} (see, e.g.,
\cite{Cow}). The term \textit{bolt of lightning} is due to Arnold \cite{Arn}%
. Marshall and O'Farrell \cite{Mar1} generalized these objects to the case
of two abstract subalgebras of the space of continuous functions defined on
a compact Hausdorff space. They gave many central properties of bolts and
functionals associated with them.

Let us now define \textit{extremal bolts}.

\bigskip

\textbf{Definition 2.2.} \textit{A finite or infinite bolt \textit{$%
\{x_{1},x_{2},...\}$} is said to be extremal for a function $f\in C(X)$ if $%
f(x_{i})=(-1)^{i}\left\Vert f\right\Vert ,i=1,2,...$ or $f(x_{i})=(-1)^{i+1}%
\left\Vert f\right\Vert ,$ $i=1,2,...$}

\bigskip

We continue with the notion of image of a finite signed measure $\mu $ and a
measure space $(U,\mathcal{A},\mu ).$ Let $F$ be a mapping from the set $U$
to the set $T.$ Then a measure space $(T,\mathcal{B},\nu )$ is called an
image of the measure space $(U,\mathcal{A},\mu )$ if the measurable sets $%
B\in \mathcal{B}$ are the subsets of $T$ such that $F^{-1}(B)\in \mathcal{A}$
and

\begin{equation*}
\nu (B)=\mu (F^{-1}(B)),\text{ for all }B\in \mathcal{B}\text{.}
\end{equation*}%
The measure $\nu $ is called an image of $\mu $ and denoted by $F\circ \mu $%
. Clearly,

\begin{equation*}
\left\Vert F\circ \mu \right\Vert \leq \left\Vert \mu \right\Vert ,
\end{equation*}%
since under mapping $F$ there is a possibility of mixing up the images of
those sets on which $\mu $ is positive with those where it is negative.
Besides, note that if a bounded function $g:T\rightarrow $ $\mathbb{R}$ is $%
F\circ \mu $-measurable, then the composite function $[g\circ
F]:U\rightarrow $ $\mathbb{R}$ is $\mu $-measurable and

\begin{equation*}
\int_{U}[g\circ F]d\mu =\int_{T}gd[F\circ \mu ].\eqno(2.1)
\end{equation*}

To prove our main result we need the following auxiliary lemmas and Singer's
theorem on characterization of a best approximation from a subspace of $C(X)$%
. By $C^{\ast }(X)$ we denote the class of regular real-valued measures of
finite total variation defined on Borel subsets of $X.$

\bigskip

\textbf{Lemma 2.1.} \textit{A measure $\mu \in C^{\ast }(X)$ is orthogonal
to $A_{1}+A_{2}$ if and only if}

\begin{equation*}
s\circ \mu \equiv 0\text{ and }p\circ \mu \equiv 0\text{.}
\end{equation*}%
\textit{That is, for any Borel subsets $E_{i}\subset X_{i}$, $i=1,2$, $\mu
(s^{-1}(E_{1}))=0$ and $\mu (p^{-1}(E_{2}))=0$.}

\bigskip

The proof of this lemma easily follows from (2.1).

\bigskip

\textbf{Lemma 2.2.} \textit{The quotient spaces $X_{1}$ and $X_{2}$ are
metrizable.}

\bigskip

This lemma is a consequence of the following two facts:

\bigskip

1) Let $A$ be a family of functions continuous on a compact space $X$ and $r$
an equivalence relation defined by $A$:

\begin{equation*}
x\overset{r}{\sim }y\text{ if }f(x)=f(y)\text{ for all }f\in A.
\end{equation*}%
Then the saturation $r(F)$ of any closed set $F\subset X$ ($r(F)\overset{def}%
{=}\bigcup\nolimits_{x\in F}r(x)$; $r(x)$ is the equivalence class of $x$)
is closed and hence the canonical projection $\pi :X\rightarrow X/r$ is a
closed mapping (see \cite[p. 54]{12}).

\bigskip

2) Let $\pi $ be a closed continuous mapping of a metric space $X$ onto a
topological space $Y$. Then the following statements are all equivalent (see
\cite{St} and \cite[Theorem 5.5]{Her}):

(a) $Y$ satisfies the first countability axiom.

(b) $\mathbf{bd}(\pi ^{-1}(y))$ is compact for each $y\in Y$ (here $\mathbf{%
bd}(\pi ^{-1}(y))$ denotes the boundary of $\pi ^{-1}(y)$).

(c) $Y$ is metrizable.

\bigskip

The following theorem plays an essential role in the proof of our main
result.

\bigskip

\textbf{Theorem 2.1.} (see Singer \cite{S}) \textit{Let $X$ be a compact
space, $M$ be a linear subspace of $C(X)$, $f\in C(X)\backslash M$ and $%
u_{0}\in M.$ Then $u_{0}$ is a best approximation in $M$ to $f$ if and only
if there exists a regular Borel measure $\mu $ on $X$ such that}

\textit{(1) The total variation $\left\Vert \mu \right\Vert =1$;}

\textit{(2) $\mu $ is orthogonal to the subspace $M$, that is, $%
\int_{X}ud\mu =0$ for all $u\in M$;}

(\textit{3) For the Jordan decomposition $\mu =\mu ^{+}-\mu ^{-}$,}
\begin{equation*}
f(x)-u_{0}(x)=\left\{
\begin{array}{c}
\left\Vert f-u_{0}\right\Vert \text{ for }x\in \text{supp}(\mu ^{+})\text{,}
\\
-\left\Vert f-u_{0}\right\Vert \text{ for }x\in \text{supp}(\mu ^{-})\text{,}%
\end{array}%
\right.
\end{equation*}%
\textit{where supp$(\mu ^{+})$ and supp$(\mu ^{-})$ are closed supports of
the positive measures $\mu ^{+}$ and $\mu ^{-}$, respectively.}

\bigskip

Our main result is the following theorem.

\bigskip

\textbf{Theorem 2.2.} \textit{Assume $X$ is a compact metric space. A
function $u_{0}\in A_{1}+A_{2}$ is a best approximation to a function $f\in
C(X)$ if and only if there exists a closed or infinite bolt extremal for the
function $f-u_{0}$.}

\bigskip

\textbf{Proof.} \textit{Necessity.} Assume $u_{0}$ is a best approximation
from $A_{1}+A_{2}$ to $f$. Since $A_{1}+A_{2}$ is a subspace of $C(X)$, we
have a regular Borel measure $\mu $ satisfying the conditions (1)-(3) of
Theorem 2.1.

Take any point $x_{0}$ in supp$(\mu ^{+})$ and\ consider the point $%
y_{0}=s(x_{0})$ in $X_{1}$. Since by Lemma 2.2, $X_{1}$ is metrizable (hence
first countable), there is a nested countable open neighborhood basis at $%
y_{0}$. Denote this basis by $\{O_{n}(y_{0})\}_{n=1}^{\infty }$. For each $n$%
, $\mu ^{+}\left[ s^{-1}(O_{n}(y_{0}))\right] >0$, since $%
s^{-1}(O_{n}(y_{0}))$ is an open set containing $x_{0}$. By Lemma 2.1, $\mu %
\left[ s^{-1}(O_{n}(y_{0}))\right] =0$. Therefore, $\mu ^{-}\left[
s^{-1}(O_{n}(y_{0}))\right] >0$. It follows that for each $n$, the
intersection $s^{-1}(O_{n}(y_{0}))\cap $supp$(\mu ^{-})$ is not empty. Take
now any points $z_{n}\in s^{-1}(O_{n}(y_{0}))\cap $supp$(\mu ^{-})$, $%
n=1,2,...$ Since supp$(\mu ^{-})$ is sequentially compact (as a closed set
in a compact metric space), the sequence $\{z_{n}\}_{n=1}^{\infty }$ or a
subsequence of it converges to a point $x_{1}$ in supp$(\mu ^{-})$. We may
assume without loss of generality that $z_{n}\rightarrow x_{1}$, as $%
n\rightarrow \infty $. Since for any $n$, $s(z_{n})\in O_{n}(y_{0})$ and $%
\{O_{n}(y_{0})\}_{n=1}^{\infty }$ is a nested neighborhood basis, we obtain
that $s(z_{n})\rightarrow y_{0}$, as $n\rightarrow \infty $. On the other
hand, since $s$ is continuous, $s(z_{n})\rightarrow s(x_{1})$, as $%
n\rightarrow \infty $. It follows that $s(x_{1})=y_{0}=s(x_{0})$. Note that $%
x_{0}\in $supp$(\mu ^{+})$ and $x_{1}\in $supp$(\mu ^{-})$.

Changing $s$ and $\mu ^{+}$ to $p$ and $\mu ^{-}$, correspondingly, repeat
the above process with the point $y_{1}=p(x_{1})$ and a nested countable
neighborhood basis at $y_{1}$. Then we obtain a point $x_{2}\in $supp$(\mu
^{+})$ such that $p(x_{2})=p(x_{1}).$ Continuing this process, we can
construct points $x_{3}\in $supp$(\mu ^{-})$ , $x_{4}\in $supp$(\mu ^{+})$ ,
and so on. Note that \ the set of all constructed points $x_{i}$, $%
i=0,1,..., $ forms a bolt. By Theorem 2.1, this bolt is extremal for the
function $f-u_{0}$.

\bigskip

\textit{Sufficiency.} The main idea in this part is the application of the
Banach-Alaoglu theorem on weak$^{\text{*}}$ sequential compactness of the
closed unit ball in $E^{\ast }$ for a separable Banach space $E$ (see, e.g.,
Rudin \cite[p. 66]{Rud}). Note that since $X$ is a compact metric space, the
space $C(X)$ is separable. Thus the closed unit ball $B$ of the continuous
dual of $C(X)$ is sequentially compact, which means that any sequence in $B$
has a convergent subsequence converging to a point in $B$.

With each bolt $l=\{x_{1},...,x_{n}\}$ with respect to $A_{1}$ and $A_{2}$,
we associate the following bolt functional

\begin{equation*}
r_{l}(F)=\frac{1}{n}\sum_{i=1}^{n}(-1)^{n+1}F(x_{i}).
\end{equation*}

It is an exercise to check that $r_{l}$ is a linear bounded functional on $%
C(X)$ with the norm $\left\Vert r_{l}\right\Vert \leq 1$ and $\left\Vert
r_{l}\right\Vert =1$ if and only if the set of points $x_{i}$ with odd
indices $i$ does not intersect with the set of points with even indices.
Besides, if $l$ is closed, then $r_{l}\in (A_{1}+A_{2})^{\perp },$ where $%
(A_{1}+A_{2})^{\perp }$ is the annihilator of the subspace $%
A_{1}+A_{2}\subset C(X).$ If $l$ is not closed, then $r_{l}$ is generally
not an annihilating functional. However, it satisfies the following
important inequality

\begin{equation*}
\left\vert r_{l}(v_{i})\right\vert \leq \frac{2}{n}\left\Vert
v_{i}\right\Vert ,\eqno(2.2)
\end{equation*}%
for all $v_{i}\in A_{i}$, $i=1,2$. This inequality means that for bolts $l$
with sufficiently large number of points, $r_{l}$ behaves like an
annihilating functional on each $A_{i}$, and hence on $A_{1}+A_{2}$. To see
the validity of (2.2) it is enough to recall that $v_{1}=$ $g\circ s$, $%
v_{2}=$ $h\circ p$ and consider the chain of equalities $%
g(s(x_{1}))=g(s(x_{2})),$ $g(s(x_{3}))=g(s(x_{4})),...$(or $%
g(s(x_{2}))=g(s(x_{3})),$ $g(s(x_{4}))=g(s(x_{5})),...$) for $%
v_{1}(x)=g(s(x))$ and similar equalities for $v_{2}(x)=h(p(x))$.

Returning to the sufficiency part of the theorem, note that there may be two
cases. The first case happens when there exists a closed bolt $%
l=\{x_{1},...,x_{2n}\}$ extremal for $f-u_{0}.$ In this case, it is not
difficult to verify that $u_{0}$ is a best approximation. Indeed, on the one
hand, the following equalities are valid:

\begin{equation*}
\left\vert r_{l}(f)\right\vert =\left\vert r_{l}(f-u_{0})\right\vert
=\left\Vert f-u_{0}\right\Vert .
\end{equation*}%
On the other hand, for any function $u\in A_{1}+A_{2}$, we have

\begin{equation*}
\left\vert r_{l}(f)\right\vert =\left\vert r_{l}(f-u)\right\vert \leq
\left\Vert f-u\right\Vert .
\end{equation*}%
Thus, $\left\Vert f-u_{0}\right\Vert \leq \left\Vert f-u\right\Vert $ for
any $u\in A_{1}+A_{2}$. That is, $u_{0}$ is a best approximation.

The second case is the existence of an infinite bolt $l=\{x_{1},x_{2},...\}$
extremal for $f-u_{0}$. In this case we proceed as follows. From $l$ we form
the finite bolts $l_{k}=\{x_{1},...,x_{k}\},$ $k=1,2,...$, and consider the
bolt functionals $r_{l_{k}}$. For the ease of notation, let us put $%
r_{k}=r_{l_{k}}.$ The sequence $\{r_{_{k}}\}_{k=1}^{\infty }$ is contained
in the closed unit ball of the dual space $C^{\ast }(X).$ By the
Banach-Alaoglu theorem, the sequence $\{r_{_{k}}\}_{k=1}^{\infty }$ must
have weak$^{\text{*}}$ cluster points. Assume $r^{\ast }$ is one of them.
Without loss of generality we may assume that $r_{k}\overset{weak^{\ast }}{%
\longrightarrow }r^{\ast },$ as $k\rightarrow \infty .$ From (2.2) it
follows that $r^{\ast }(v_{1}+v_{2})=0,$ for any $v_{i}\in A_{i},$ $i=1,2$.
That is, $r^{\ast }$ belongs to the annihilator of the subspace $A_{1}+A_{2}$%
. Since we have also $\left\Vert r^{\ast }\right\Vert \leq 1,$ it follows
that

\begin{equation*}
\left\vert r^{\ast }(f)\right\vert =\left\vert r^{\ast }(f-u)\right\vert
\leq \left\Vert f-u\right\Vert ,\eqno(2.3)
\end{equation*}%
for all functions $u\in A_{1}+A_{2}.$ On the other hand, since the infinite
bolt $\{x_{1},x_{2},...\}$ is extremal for $f-u_{0},$

\begin{equation*}
\left\vert r_{k}(f-u_{0})\right\vert =\left\Vert f-u_{0}\right\Vert ,\text{ }%
k=1,2,...
\end{equation*}%
Hence

\begin{equation*}
\left\vert r^{\ast }(f)\right\vert =\left\vert r^{\ast }(f-u_{0})\right\vert
=\left\Vert f-u_{0}\right\Vert .\eqno(2.4)
\end{equation*}%
From (2.3) and (2.4) we obtain that

\begin{equation*}
\left\Vert f-u_{0}\right\Vert \leq \left\Vert f-u\right\Vert ,
\end{equation*}%
for all $u\in A_{1}+A_{2}.$ This means that $u_{0}$ is a best approximation
to $f$. Q.E.D.

\bigskip

\textbf{Remark 1.} In \cite{AI}, Theorem 2.2 was proved under additional
assumption that the algebras have the $C$-property, that is, for any $w\in
C(X)$, the functions

\begin{eqnarray*}
g_{1}(a) &=&\max_{\substack{ x\in X  \\ s(x)=a}}w(x),\text{ }g_{2}(a)=\min
_{\substack{ x\in X  \\ s(x)=a}}w(x),\text{ }a\in X_{1}, \\
h_{1}(b) &=&\max_{\substack{ x\in X  \\ p(x)=b}}w(x),\text{ }h_{2}(b)=\min
_{\substack{ x\in X  \\ p(x)=b}}w(x),\text{ }b\in X_{2}
\end{eqnarray*}%
are continuous.

\bigskip

\textbf{Remark 2.} Note that in the special case when $X\subset \mathbb{R}%
^{2}$ and $s,p$ are the coordinate functions, a Chebyshev type alternation
theorem was first obtained by Khavinson \cite{Kh}. In \cite{7}, similar
alternation theorems were proved for ridge functions and certain function
compositions.

\bigskip

\textbf{Remark 3.} Note that characterization of a best approximation from a
sum of more than two subalgebras $A_{1},...,A_{k}$ of $C(X)$ seems to be
beyond the scope of the methods discussed herein. A bolt with respect to two
algebras $A_{1}$ and $A_{2}$ is constructed as a sequence of points $%
\{x_{1},x_{2},...\}$ with the links $x_{i}x_{i+1}$ traveling alternatively
in equivalence classes of the relations $R_{1}$ and $R_{2}$ (see above). In
this case, the bolt functional $r_{l}$ has important property (2.2), which
leads to the functional $r^{\ast }$ annihilating all elements of the sum $%
A_{1}+A_{2}$. The problem becomes complicated when the number of summands in
the sum $A_{1}+\cdot \cdot \cdot +A_{k}$ is more than two. The simple
generalization of bolts demands a sequence of points $\{x_{1},x_{2},...\}$
with the links $x_{i}x_{i+1}$ traveling in three or more alternating
equivalence classes. But in this case, the number $2$ in (2.2) grows
unboundedly as $n$ tends to infinity and we cannot arrive at any
annihilating functional like $r^{\ast }$. For $k\geq 3$ we do not know a
reasonable description of a sequence of points $\{x_{1},x_{2},...\}$ and
functionals $r_{l_{n}}$, associated with the first $n$ points $%
x_{1},...,x_{n}$, such that any weak$^{\text{*}}$ cluster point of the
sequence $\{r_{l_{n}}\}_{n=1}^{\infty }$ is orthogonal to the sum $%
A_{1}+\cdot \cdot \cdot +A_{k}$. We refer the interested reader to Sternfeld \cite{Ste}
for discussions on differences between the cases of two and more than two
algebras.

\bigskip

\begin{center}
{\large \textbf{3. Competing interests declaration}}
\end{center}

The authors declare none.


\bigskip

\begin{thebibliography}{99}
\bibitem{Arn} V. I. Arnold, On functions of three variables, (Russian)%
\textit{\ Dokl. Akad. Nauk SSSR} \textbf{114} (1957), 679-681; English
transl. in \textit{Amer. Math. Soc. Transl. }\textbf{28} (1963), 51-54.

\bibitem{AI} A. Kh. Asgarova, V. E. Ismailov, A Chebyshev-type theorem
characterizing best approximation of a continuous function by elements of
the sum of two algebras, (Russian) \textit{Mat. Zametki} \textbf{109}
(2021), no. 1, 19-26.

\bibitem{Bro} B. Brosowski, A. R. da Silva, A general alternation theorem,
Approximation theory (Memphis, TN, 1991), 137-150, \textit{Lecture Notes in
Pure and Appl. Math.}, 138, Dekker, New York, 1992.

\bibitem{Buck} R. C. Buck, Alternation theorems for functions of several
variables, \textit{J. Approx. Theory} \textbf{1} (1968), 325-334.

\bibitem{Cow} R. C. Cowsik, A. Klopotowski, M. G. Nadkarni, When is $%
f(x,y)=u(x)+v(y)$?, \textit{Proc. Indian Acad. Sci. Math. Sci.} \textbf{109}
(1999), 57-64.

\bibitem{4} S. P. Diliberto, E. G. Straus, On the approximation of a
function of several variables by the sum of functions of fewer variables,
\textit{Pacific J. Math.} \textbf{1} (1951), 195-210.

\bibitem{Her} R. A. Herman, \textit{Quotients of metric spaces}, A.B.,
Grinnell College, 1966.

\bibitem{7} V. E. Ismailov, \textit{Ridge functions and applications in
neural networks}, Mathematical Surveys and Monographs 263, American
Mathematical Society, 186 pp.

\bibitem{12} S. Ya. Khavinson, \textit{Best approximation by linear
superpositions (approximate nomography),} Translated from the Russian
manuscript by D. Khavinson. Translations of Mathematical Monographs, 159.
American Mathematical Society, Providence, RI, 1997, 175 pp.

\bibitem{Kh} S. Ya. Khavinson (S. Ja. Havinson), A Chebyshev theorem for the
approximation of a function of two variables by sums of the type $\varphi
\left( {x}\right) +\psi \left( {y}\right) ,$ \textit{Izv. Acad. Nauk. SSSR
Ser. Mat.} \textbf{\ 33} (1969), 650-666; English tarnsl. in \textit{Math.
USSR Izv.} \textbf{3} (1969), 617-632.

\bibitem{13} W. A. Light, E. W. Cheney, On the approximation of a bivariate
function by the sum of univariate functions. \textit{J. Approx. Theory}
\textbf{29} (1980), 305--322.

\bibitem{Mar1} D. E. Marshall, A. G. O'Farrell, Approximation by a sum of
two algebras. The lightning bolt principle, \textit{J. Funct. Anal. }\textbf{%
52} (1983), 353-368.

\bibitem{Mar2} D. E. Marshall, A. G. O'Farrell, Uniform approximation by
real functions, \textit{\ Fund. Math.} \textbf{104} (1979), 203-211.

\bibitem{Nat} I. P. Natanson, \textit{Constructive function theory, Vol. I.
Uniform approximation,} Translated from the Russian by Alexis N. Obolensky
Frederick Ungar Publishing Co., New York 1964, 232 pp.

\bibitem{Rud} W. Rudin, \textit{Functional analysis}, McGraw-Hill Series in
Higher Mathematics. McGraw-Hill Book Co., 1973, 397 pp.

\bibitem{S} I. Singer, \textit{The theory of best approximation and
functional analysis}. Conference Board of the Mathematical Sciences Regional
Conference Series in Applied Mathematics, No. 13. Society for Industrial and
Applied Mathematics, Philadelphia, Pa., 1974, 95 pp.

\bibitem{Ste} Y. Sternfeld, Uniform separation of points and measures and
representation by sums of algebras, \textit{Israel J. Math.} \textbf{55}
(1986), no. 3, 350--362.

\bibitem{St} A. H. Stone, Metrizability of decomposition spaces, \textit{%
Proc. Amer. Math. Soc.} \textbf{7} (1956), 690-700.

\bibitem{St1} M. H. Stone, Applications of the Theory of Boolean Rings to
General Topology, \textit{Trans. Amer. Math. Soc.} \textbf{41} (1937),
375--481.

\bibitem{St2} M. H. Stone, The generalized Weierstrass approximation
theorem, \textit{Math. Mag.} \textbf{21} (1948), 167--184, 237--254.
\end{thebibliography}
\end{document}